\documentclass{article}

\usepackage{arxiv}

\usepackage[utf8]{inputenc}
\usepackage[T1]{fontenc}
\usepackage{hyperref}
\usepackage{url}
\usepackage{booktabs}
\usepackage{amsfonts}
\usepackage{nicefrac}
\usepackage{microtype}
\usepackage{graphicx}

\usepackage[numbers,sort&compress]{natbib}
\usepackage{doi}

\usepackage{amsmath,mathtools,amssymb}
\usepackage{xcolor}
\usepackage{enumitem}
\usepackage{algorithm}
\usepackage{multicol}

\newtheorem{assumption}{\bf{Assumption}}

\newtheorem{remark}{\bf{Remark}}

\newcommand{\xddots}{%
  \raise 4pt \hbox {.}
  \mkern 6mu
  \raise 1pt \hbox {.}
  \mkern 6mu
  \raise -2pt \hbox {.}
}

\DeclareUnicodeCharacter{221E}{\ensuremath{\infty}}

\graphicspath{{images/}}

\title{Generalized String-Stability Criteria for Consensus Protocols}

\author{{Sridhar Babu Mudhangulla} \\
	Department of Electrical and Computer Engineering\\
	Florida State University\\
	Tallahassee, FL 32303 \\
	\texttt{sm19ch@fsu.edu} \\
	\And
	{Olugbenga Moses Anubi} \\
	Department of Electrical and Computer Engineering\\
	Florida State University\\
	Tallahassee, FL 32303 \\
	\texttt{oanubi@fsu.edu} \\
}

\date{} 

\renewcommand{\headeright}{}
\renewcommand{\undertitle}{}
\renewcommand{\shorttitle}{}

\hypersetup{
pdftitle={Generalized String-Stability Criteria for Consensus Protocols},
pdfsubject={Control Theory, Multi-Agent Systems, string stability, consensus},
pdfauthor={Sridhar Babu Mudhangulla, Olugbenga Moses Anubi},
pdfkeywords={Multi-agent systems, Consensus algorithms, Leader-follower networks, string stability, Input-to-error stability, Distributed control},
}

\begin{document}
\maketitle

\begin{abstract}
This paper presents a unified string-stability framework for leader--follower multi-agent systems governed by first-, second-, and $m^{\text{th}}$-order consensus protocols operating under an $r$-predecessor directed communication topology. While string stability has been extensively studied for specific vehicle models and individual consensus protocols, existing results remain fragmented across protocol orders and do not identify the fundamental factors governing disturbance amplification or attenuation. This work shows that, for all consensus orders, string stability is dictated solely by the communication richness $r$, while the protocol order $m$ influences the transient and frequency-domain response, including damping and mid- to high-frequency behavior. In particular, the low-frequency gain of the disturbance propagation coefficient is inversely proportional to $r$ for every $m$, implying that higher-order consensus cannot alter the structural low-frequency limitation imposed by insufficient communication, although it significantly influences transient behavior depending on communication richness and gain selection. Under the adopted $\mathcal{H}_\infty$-based string-stability definition and the present framework, string stability is achievable if and only if $r \ge 2$. This establishes a structural–dynamic separation principle that unifies and generalizes classical platoon results, providing new insight into the interplay between topology and controller design in cooperative driving and multi-agent coordination. The framework is developed under idealized identical-agent and fixed-topology assumptions, providing a baseline for future robust extensions. Numerical simulations corroborate the analysis and further illustrate that the benefits of higher-order consensus become effective only when sufficient predecessor information is available.
\end{abstract}

\keywords{Multi-agent systems, consensus algorithms, leader--follower networks, string stability, distributed control.}

\section{Introduction}
Cooperative control of networked systems has emerged as a cornerstone of modern engineering, enabling autonomous agents to coordinate using only local interactions. Applications span connected and automated vehicles (CAVs), robotic formations, sensor networks, and large-scale cyber–physical systems (see, e.g., \cite{olfati2004consensus,ren2007information}). In many such systems, agents are arranged along a spatial axis, requiring the collective response to upstream disturbances to remain stable and well-behaved. In vehicle platoons and robotic convoys in particular, this requirement is formalized through string stability, which requires that disturbances introduced at the leader or at any agent do not amplify as they propagate downstream. Violations of string stability lead to unsafe spacing fluctuations, increased energy consumption, degraded comfort, and ultimately limit the scalability of cooperative systems \cite{swaroop2002string}.

Classical adaptive cruise control (ACC), which relies solely on local measurements, was shown to be inherently string-unstable \cite{barooah2005error}, motivating cooperative adaptive cruise control (CACC) architectures that exploit communication among vehicles \cite{ploeg2013lp}. Numerous studies have examined string stability in CACC-equipped platoons \cite{vehicular_stdli_2017, performance_hansson_2024, overview_salek_2025}, revealing a strong interplay between communication topology and control design. These results raise a fundamental question central to distributed coordination: How do local feedback laws and information flow jointly determine disturbance propagation in networked systems?

Consensus protocols offer a natural modeling framework for distributed coordination. First-order consensus regulates position disagreement, second-order consensus incorporates velocity alignment, and higher-order extensions generalize to more complex agent dynamics \cite{ren2008consensus,wieland2008consensus,overview_knorn_2016}. While structurally similar to cooperative longitudinal control, consensus theory traditionally focuses on asymptotic agreement and does not explicitly address disturbance propagation across spatially ordered networks \cite{constrained_zhou_2019, fully_mazouchi_2021, output_zhu_2025}. Consequently, consensus theory and string-stability analysis have progressed largely independently.

In the platooning literature, string stability is typically analyzed in the frequency domain for specific vehicle models and fixed communication structures \cite{besselink2017string,pare2019networked,disturbance_seiler_2004}. Conversely, recent consensus-based studies \cite{distributed_huang_2020, multiagent_jiang_2022} emphasize spectral convergence properties and graph connectivity but rarely characterize how disturbances evolve along an $r$-predecessor communication topology. As a result, existing string-stability results remain fragmented across controller types and do not reveal whether disturbance attenuation is governed primarily by controller order, communication richness, or their interaction. This gap motivates the need for a unified treatment that connects arbitrary-order consensus protocols with structured communication topologies. Recent developments have further extended consensus and cooperative control toward robust, adaptive, and delay-tolerant architectures, including distributed output-feedback control under delays, consensus with heterogeneous time-varying input and communication delays, and output-consensus problems for heterogeneous nonlinear networked systems \cite{distributed_huang_2020,multiagent_jiang_2022,output_zhu_2025}. These directions emphasize practical scalability and robustness, complementing the present work’s focus on fundamental disturbance-propagation limits under structured communication topologies.

This paper addresses this gap by developing a unified string-stability framework for leader--follower multi-agent systems governed by first-, second-, and general $m^{\text{th}}$-order consensus protocols operating under an $r$-predecessor directed topology. We leverage a frequency-domain characterization of disturbance propagation to show that the low-frequency gain of the closed-loop mapping from the leader's disturbance to any follower equals $1/r$ for all consensus orders. This reveals a fundamental structural limitation: under the adopted $\mathcal{H}_\infty$-based string-stability definition and the present framework, string stability is achievable if and only if $r \ge 2$, independent of consensus order and gain tuning at low frequencies. Higher-order consensus modifies the transient and frequency-domain behavior, influencing damping, oscillations, and mid- to high-frequency attenuation, while the structural low-frequency limit remains unchanged. This result unifies and generalizes insights from classical platooning and consensus control, clarifying the separation between structural and dynamic disturbance attenuation mechanisms. Moreover, the simulations reveal that the benefits of higher-order consensus become effective only when sufficient predecessor information is available, highlighting an interaction between communication richness and dynamic shaping.

The main contributions of this paper are as follows:
\begin{enumerate}
\item A unified frequency-domain framework for analyzing string stability of leader--follower multi-agent systems under first-, second-, and general $m^{\text{th}}$-order consensus protocols with $r$-predecessor communication.
\item A structural string-stability condition proving that the low-frequency disturbance gain is $1/r$ for all consensus orders, implying string stability if and only if $r \ge 2$.
\item A characterization of the role of consensus order, showing that increasing $m$ introduces additional dynamic shaping that influences damping, oscillations, and mid- to high-frequency attenuation, without affecting the low-frequency structural limit.
\item A conceptual unification of consensus theory and classical string-stability results, clarifying the interplay between communication richness and dynamic protocol order.
\item Numerical validation illustrating how $r$ and $m$ jointly influence disturbance propagation, damping, and transient behavior.
\end{enumerate}

\section{Notations and Preliminaries}
\subsection{Notations}
The set of real numbers is denoted by $\mathbb{R}$, and $\mathbb{R}^n$ and $\mathbb{R}^{n\times m}$ denote real vectors of dimension $n$ and real matrices of size $n\times m$, respectively. $\mathbb{R}_+$ denotes the set of positive real numbers. The index set $\mathcal{I}_n=\{1,2,\ldots,n\}$. Boldface lower-case symbols (e.g., $\mathbf{x}$) denote vectors, and normal upper-case symbols (e.g., $X$) denote matrices. $x_i^{(k)}$ denotes the $k^{\text{th}}$ derivative of $x_i$ with $x_i^{(0)} = x_i$. $\mathcal{L}_2=\left\{x: \mathbb{R}_+ \rightarrow \mathbb{R} \mid \int_0^{\infty} x^2(t)\,dt<\infty\right\}$ denotes the space of all square-integrable signals. Given $x \in \mathcal{L}_2, X(s) \triangleq \int_0^{\infty} e^{-s t} x(t)\,dt$ denotes its Laplace transform. Given a system transfer function $G(j \omega)$, its induced $\mathcal{H}_\infty$ norm is given by \cite{zhou1998essentials}
\begin{align*}
\|G\|_{\mathcal{H}_\infty}=\sup_{\omega\ge0}|G(j \omega)|.
\end{align*}

\subsection{Graph Theory}
The tuple $\mathcal{G}=(\mathcal{V},\mathcal{E})$ denotes a directed graph with vertex set $\mathcal{V}=\{v_1,\ldots,v_n\}$ and edge set $\mathcal{E}\subseteq\mathcal{V}\times\mathcal{V}$. An edge $(v_j,v_i)\in\mathcal{E}$ indicates that agent $i$ receives information from agent $j$. 

The in-neighbor set of agent $i$ is $\mathcal{N}_i=\left\{j:\left(v_j, v_i\right) \in \mathcal{E}\right\}$, and the in-degree of $v_i$ is $d_i=|\mathcal{N}_i|$. The graph Laplacian is given by $L = D - A$, where $D=\operatorname{diag}(d_1,\ldots,d_n)$ is the degree matrix and $A=[a_{ij}]$ with
\begin{align*}
a_{i j}= \begin{cases}1, & \left(v_j, v_i\right) \in \mathcal{E}, \\ 0, & \text { otherwise },\end{cases}
\end{align*}
is the adjacency matrix. A spanning tree is a directed tree that contains all the nodes and some edges. The digraph $\mathcal{G}$ contains a spanning tree if the corresponding Laplacian matrix has a simple zero eigenvalue associated with the right eigenvector $\mathbf{1}$ (\cite{zegers2016consensus,godsil2013algebraic}).

\subsection{String Stability}
\begin{figure} 
\centerline{\includegraphics[width = 0.48\textwidth]{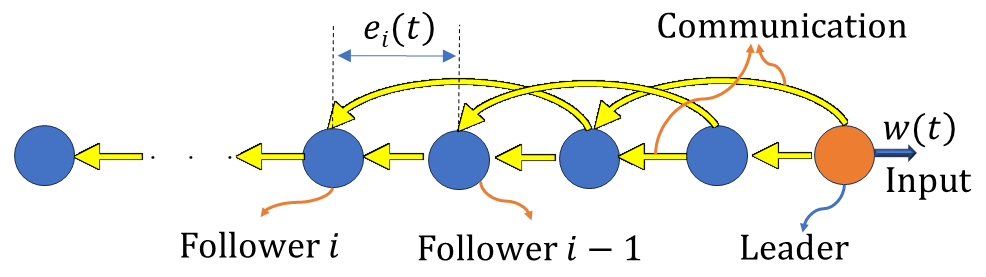}}
\caption{Directed communication topology of a platoon. The leader (orange) injects the disturbance, followers (blue) receive information through the incoming yellow links.}
\label{fig:platoon}
\end{figure}

In leader--follower formations (see Figure~\ref{fig:platoon}), string stability characterizes whether disturbances applied at the leader amplify as they propagate along the formation. Under the assumption that the leader’s disturbance acts as an external input and all followers start from zero initial conditions \cite{feng2019string}, string stability is naturally defined using the induced $\mathcal{H}_{\infty}$ norm of the disturbance-to-output transfer function.

Let $G_{0,i}(s)$ denote the transfer function from the leader (index $0$) to follower $i$. A platoon is string stable \cite{naus2010string} if
\begin{align*}
\max_{i\in\mathcal{I}_n}\left\|G_{0,i}\right\|_{\mathcal{H}_\infty}\leq1,
\end{align*}
that is, the leader-induced disturbance is not amplified at any follower. This definition is appropriate when the leader's input is treated as an external excitation, and all vehicles begin from zero initial conditions. The remainder of this paper adopts the $\mathcal{H}_{\infty}$-based notion of string stability. For alternative definitions and a comprehensive taxonomy of string stability, see \cite{feng2019string}. In subsequent analysis, the leader’s influence is modeled as an equivalent disturbance input to the follower error dynamics.

\section{Problem Setup} \label{sec:model}
Consensus algorithms describe distributed control laws that enable a network of agents to asymptotically agree using only local interactions. For agents modeled as chains of $m$ integrators, consensus protocols regulate relative derivatives of order $0$ through $m-1$. Here we develop a unified modeling framework for first-, second-, and general $m^{\text {th }}$-order consensus protocols and derive the closed-loop dynamics under the $r$-predecessor communication structure that forms the basis for subsequent string stability analysis.

The present framework adopts several assumptions to enable closed-form analysis:
\begin{assumption} \label{ass:dynamics}
    Identical follower dynamics modeled as chains of integrators.
\end{assumption}
\begin{assumption} \label{ass:topology}
    Fixed directed $r$-predecessor information topology.
\end{assumption}
\begin{assumption} \label{ass:communication}
    Ideal communication without delays, dropout, or quantization.
\end{assumption}
\begin{assumption} \label{ass:gains}
    Homogeneous protocol gains across agents.
\end{assumption}
These assumptions isolate the structural relationship between communication richness and disturbance propagation. Relaxing these assumptions constitutes an important direction for future work.

\subsection{General $m^{\text{th}}$-order Consensus Protocol}
Considering a multi-agent system consisting of $n$ agents, each with $m$ cascaded integrator dynamics of the form $x_{i}^{(m)} = u_i$, $i \in \mathcal{I}_n$, which can be written in first-order form as
\begin{align}  \label{equ:agent_dynamics_morder}
    \dot{x}_i^{(0)} &= x_i^{(1)} \nonumber \\
    &\vdots \nonumber \\
    \dot{x}_i^{(m-2)} &= x_i^{(m-1)} \nonumber \\
    \dot{x}_i^{(m-1)} &= u_i, \text{ for all } i \in \mathcal{I}_n.
\end{align}
The $m^{\text{th}}$-order consensus protocol \cite{wieland2008consensus,ren2006high} is given as
\begin{align}\label{equ:nth_protocol}
    u_i &= -\sum_{j=1}^n a_{ij} \left[\sum_{k=0}^{m-1} \gamma_k \big(x_{i}^{(k)} - x_{j}^{(k)}\big)\right],
\end{align}
where $x_i^{(k)}\in\mathbb{R}$, $k=0,\dots,m-1$, denote the states of the $i^{\text{th}}$ agent, $i \in \mathcal{I}_n$, $a_{ij}$ are adjacency weights, and $\gamma_k > 0, k=0,\dots,m-1$, are protocol gains, which are independent of $i$ due to Assumption~\ref{ass:gains}.

Let $\mathbf{x}^{(k)} = [{x}_1^{(k)}, \dots, {x}_n^{(k)}]^\top$ and $\mathbf{u} = [{u}_1, \dots, {u}_n]^\top$ then algorithm \eqref{equ:nth_protocol} for $n$ agents can be written in vector form, 
\begin{align}\label{equ:nagents_nthprotocol}
    \mathbf{u} &= -\sum_{k=0}^{m-1} \gamma_{k} L \mathbf{x}^{(k)}.
\end{align}
Furthermore, define the aggregated state vector
\begin{align*}
    \boldsymbol{\xi} = [{\mathbf{x}}^{(0)} \dots {\mathbf{x}}^{(m-1)}] ^\top \in \mathbb{R}^{mn}.
\end{align*}
Then the resulting closed-loop dynamics is given by
\begin{align} \label{equ:nordersys_mthprotocol}
    \dot{\boldsymbol{\xi}} = \mathbf{M}\boldsymbol{\xi},
\end{align}
where 
\begin{align*} 
\mathbf{M} =\begin{bmatrix}
    0_{n} &I_{n} &0_{n} &\dots &0_{n} \\
    0_{n} &0_{n} &I_{n} &\dots &0_{n} \\
    \vdots &\vdots &\vdots &\vdots &\vdots \\
    - \gamma_0 L &- \gamma_1 L &-\gamma_2 L &\dots &-\gamma_{m-1} L \end{bmatrix} \in \mathbb{R}^{mn \times mn}.
\end{align*}
The characteristic equation is obtained, using the Schur complement, as 
\begin{align*}
    \det(\lambda I_{mn} - \mathbf{M}) = 
    \prod_{i=1}^{n}\Big(\lambda^m-\sum_{k=0}^{m-1}\gamma_k \mu_i \lambda^{k}\Big) = 0,
\end{align*}
where $\mu_i$ is the $i^{\text{th}}$ eigenvalue of $-L$, $i=1,\dots,n$. Under the directed spanning tree assumption, $-L$ has one zero eigenvalue $\mu_1=0$ with right eigenvector $\mathbf{1}$, and all other eigenvalues $\mu_2, \ldots, \mu_n$ satisfy $\operatorname{Re}\left(\mu_i\right)<0$ (note that $\mu_i$ is the eigenvalue of $-L$). From the characteristic equation above, it is seen that each mode satisfies the scalar $m^{\text{th}}$ degree polynomial equation:
\begin{align} \label{eq:mode_poly}
    \lambda^m = \gamma_{m-1} \mu_i \lambda^{m-1}+\cdots+\gamma_1 \mu_i \lambda+\gamma_0 \mu_i.
\end{align}
Thus the eigenvalues of the full matrix $M$ are precisely the roots of \eqref{eq:mode_poly} for each Laplacian eigenvalue $\mu_i$. 

The network is said to achieve internal consensus if
\begin{align*}
\lim_{t\to\infty}\left|x_i^{(k)}(t)-x_j^{(k)}(t)\right|=0,
\quad
\forall i,j\in\mathcal I_n,\;\forall k=0,\ldots,m-1.
\end{align*}
Thus, in summary, the $m^{\text{th}}$-order consensus for \eqref{equ:nordersys_mthprotocol} is achieved if:
\begin{enumerate}
    \item The directed graph described by $L$ has a directed spanning tree, and
    \item The matrix $\mathbf{M}$ in \eqref{equ:nordersys_mthprotocol} has $m$ zero eigenvalues with geometric multiplicity one, and all the other eigenvalues have negative real parts (This can further be ascertained using Routh-Hurwitz or equivalent stability criteria).
\end{enumerate}

\begin{remark}
    For first-order consensus protocol, setting $m=1$, the closed loop dynamics 
    \begin{align*}
    \dot{x}_i=-\gamma_0 \sum_{j \in \mathcal{N}_i} a_{i j}\left(x_i-x_j\right),
    \end{align*}
    satisfies $\lambda=\gamma_0 \mu_i$. Thus internal stability is guaranteed if $\gamma_0>0$ and the graph has a spanning tree. Similarly, for the second-order consensus protocol ($m=2$), the closed loop dynamics 
    \begin{align*}
    \ddot{x}_i=-\gamma_1 \sum_{j \in \mathcal{N}_i} a_{i j}\left(\dot{x}_i-\dot{x}_j\right)-\gamma_0 \sum_{j \in \mathcal{N}_i} a_{i j}\left(x_i-x_j\right),
    \end{align*}
    satisfies the characteristic equation $\lambda^2-\gamma_1 \mu_i \lambda-\gamma_0 \mu_i=0$, for which the standard second-order stability condition requires positive proportional and derivative gains, i.e., $\gamma_0>0$ and $\gamma_1>0$, under the nominal consensus setting.
\end{remark}

In the remainder of the paper, we specialize the general consensus protocol to the $r$-predecessor topology, which enables a structured representation of the closed-loop dynamics and admits tractable disturbance-propagation analysis.

\subsection{Closed-Loop Error Dynamics Under $r$-Predecessor}
For the string-stability analysis, we focus on interior followers in a long leader--follower chain, for which each follower has access to its $r$ immediate predecessors:
\begin{align*}
\mathcal N_i=\{i-1,i-2,\ldots,i-r\}, \qquad i>r.
\end{align*}
Under this interior approximation, the induced communication matrix is lower-triangular, Toeplitz, and shift-invariant, with
\begin{align*}
L_{ij}=
\begin{cases}
r, & i=j,\\
-1, & j\in\mathcal N_i,\\
0, & \text{otherwise}.
\end{cases}
\end{align*}
Unlike general directed graphs, this structured topology permits exact closed-form frequency-domain expressions and a scalar disturbance-propagation model that avoids explicit eigenvalue decompositions.

Let $x_0^{(k)}(t)$, $k=0,1,\ldots,m-1$, denote the leader trajectory and its derivatives, and define the relative tracking errors
\begin{align*}
\epsilon_i^{(k)}(t)=x_i^{(k)}(t)-x_0^{(k)}(t), \qquad k=0,1,\ldots,m-1,
\end{align*}
with $\epsilon_i(t)=\epsilon_i^{(0)}(t)=x_i(t)-x_0(t)$. The $m^{\text{th}}$-order consensus protocol specialized to the $r$-predecessor structure is given by
\begin{align}
u_i(t)=-a\sum_{j=1}^{r}\sum_{k=0}^{m-1}\gamma_k\left(x_i^{(k)}(t)-x_{i-j}^{(k)}(t)\right),
\label{eq:m_protocol}
\end{align}
where $a>0$ denotes the uniform coupling strength associated with predecessor interactions. Substituting~\eqref{eq:m_protocol} into the follower dynamics yields the error system
\begin{align} \label{eq:m_error}
    \epsilon_i^{(m)}(t)=-a \sum_{j=1}^r \sum_{k=0}^{m-1} \gamma_k\left(\epsilon_i^{(k)}(t)-\epsilon_{i-j}^{(k)}(t)\right)+w(t)
\end{align}
where the disturbance is defined as $w(t) = -x^{(m)}_0(t)$. Thus, the leader motion enters the follower error dynamics as an explicit external disturbance input. Taking the Laplace transform of \eqref{eq:m_error} under zero initial conditions gives
\begin{align} \label{eq:m_error_lap}
    s^m E_i(s)=-a \sum_{j=1}^r \sum_{k=0}^{m-1} \gamma_k s^k\left(E_i(s)-E_{i-j}(s)\right)+W(s).
\end{align}
Define the consensus shaping polynomial
\begin{align} \label{eq:Qs}
    Q_m(s)=\sum_{k=0}^{m-1} \gamma_k s^k.
\end{align}
Then \eqref{eq:m_error_lap} can be rewritten as
\begin{align*} 
    \left(s^m+a r Q_m(s)\right) E_i(s)-a Q_m(s) \sum_{j=1}^r E_{i-j}(s)=W(s).
\end{align*}
or equivalently,
\begin{align} \label{eq:Ei_general}
    E_i(s)=G_m(s)\left(W(s)+a Q_m(s) \sum_{j=1}^r E_{i-j}(s)\right),
\end{align}
where 
\begin{align}
    G_m(s)=\frac{1}{s^m+ a r Q_m(s)}.
\end{align}
Defining the propagation coefficient
\begin{align}
\Phi_m(s)=a Q_m(s) G_m(s)=\frac{a Q_m(s)}{s^m+a r Q_m(s)}.
\label{eq:phi_def}
\end{align}
Then \eqref{eq:Ei_general} becomes
\begin{align}
E_i(s)=G_m(s)W(s)+\Phi_m(s)\sum_{j=1}^r E_{i-j}(s).
\label{eq:prop_coefficient}
\end{align}
The scalar transfer function $\Phi_m(s)$ governs how disturbances propagate through the formation. A related transfer function was reported in \cite{bian2019reducing} for the $r$-predecessor topology under a constant time-headway policy. In that setting, the propagation factor depends on the predecessor position because of the headway term, whereas \eqref{eq:phi_def} is independent of predecessor index and depends only on the total number of predecessors $r$. This common propagation coefficient captures the local mechanism by which disturbances are transmitted under the Toeplitz $r$-predecessor topology. Its invariance with respect to the follower index $i$ follows directly from the shift-invariant communication structure. 

Consequently, string stability can be analyzed through the frequency response magnitude $\left|\Phi_m(j\omega)\right|$. Evaluating \eqref{eq:phi_def} at zero frequency yields $\left|\Phi_m(0)\right|=\frac{1}{r}$, which is independent of the consensus order $m$ and the gains $\gamma_k$. This fundamental property underlies the unified string-stability results developed next. In contrast, the order $m$ and the coefficients $\gamma_k$ influence the transient and frequency-domain behavior, including damping, oscillations, and mid- to high-frequency attenuation, without altering the structural low-frequency limit. Moreover, the effectiveness of higher-order consensus in improving transient performance depends on sufficient communication richness, as limited predecessor information may lead to oscillatory or poorly damped responses.

\section{String Stability of Consensus Protocols} \label{sec:string_stability}
This section develops frequency-domain conditions for disturbance attenuation in leader--follower formations governed by $m^{\text{th}}$-order consensus protocols under the $r$-predecessor topology. Using the closed-loop model derived above, we show that the coefficient $\Phi_m(s)$ captures the fundamental local mechanism by which disturbances are transmitted from upstream agents to downstream followers, while the communication parameter $r$, rather than the consensus order $m$, determines the principal structural limits of attenuation.

Classically, string stability is defined in terms of the leader-to-follower transfer functions $G_{0,i}(s)$, requiring that disturbances introduced at the leader do not amplify downstream. In the structured $r$-predecessor topology considered here, these global transfer maps are recursively generated through the common local propagation coefficient $\Phi_m(s)$. Motivated by this structure, we adopt
\begin{align*}
\sup_{\omega\ge0}\left|\Phi_m(j\omega)\right|\le1
\end{align*}
as an operative propagation-based non-amplification criterion. This condition characterizes whether disturbance components grow from one stage of the formation to the next, and therefore provides a direct indicator of downstream disturbance attenuation in the present framework.

It should be noted that the above criterion is based on local propagation behavior and is not, by itself, identical to the classical leader-to-follower $\mathcal H_\infty$ condition. Nevertheless, it captures the dominant structural mechanism governing amplification in long homogeneous formations, and therefore serves as the central analytical tool for characterizing local disturbance propagation. While the supremum condition is compact, it masks important distinctions across frequency ranges. Low-frequency behavior governs spacing drift, mid-frequency behavior determines oscillatory transients, and high-frequency behavior reflects sensitivity to measurement noise. These distinctions are critical because the consensus order $m$ and communication richness $r$ influence different spectral regions in fundamentally different ways.

In contrast to the existing platooning literature, which analyzes specific vehicle models or controller types, our approach reveals how protocol order and communication richness jointly shape attenuation through the single scalar function
\begin{align*}
\Phi_m(s)=\frac{a Q_m(s)}{s^m+a r Q_m(s)}, \text{ where } Q_m(s)=\sum_{k=0}^{m-1} \gamma_k s^k.
\end{align*}
The Toeplitz structure of the $r$-predecessor Laplacian ensures that this propagation coefficient is identical for all followers, making string stability analysis transparent and explicit. We now characterize first-, second-, and general $m^{\text{th}}$-order protocols within a unified framework.

\subsection{String Stability of First-Order Consensus}
First-order consensus is the simplest and most widely used distributed control law. For $m=1$, the general consensus protocol reduces to
\begin{align*}
u_i(t)=-a \gamma_0 \sum_{j=1}^r\left(x_i(t)-x_{i-j}(t)\right),
\end{align*}
resulting in the propagation coefficient
\begin{align*}
    \Phi_1(s)=\frac{a \gamma_0}{s+a r \gamma_0},
\end{align*}
which implies that
\begin{align} \label{eq:Phi_frequency}
\left|\Phi_1(j \omega)\right|=\frac{a \gamma_0}{\sqrt{\left(a r \gamma_0\right)^2+\omega^2}}.
\end{align}
It follows from \eqref{eq:Phi_frequency} that
\begin{align*}
\left|\Phi_1(0)\right|=\frac{1}{r},
\qquad
\left|\Phi_1(j\omega)\right|\sim \frac{a\gamma_0}{|\omega|}\to0
\quad \text{as } \omega\to\infty.
\end{align*}
Thus, low-frequency disturbance attenuation depends entirely on communication richness. For $r=1$, the zero-frequency gain equals unity, and disturbances do not decay spatially; for $r\ge2$, the formation is structurally string stable. High-frequency attenuation scales as $1/\omega$, and no resonant peak arises due to the first-order closed-loop dynamics. These observations recover the classical view that first-order consensus provides no disturbance attenuation under single-predecessor communication, but becomes string stable once each follower receives information from at least two predecessors.

\subsection{String Stability of Second-Order Consensus}
Second-order consensus introduces velocity alignment. For each follower $i$, the control law is
\begin{align*}
u_i(t)=-a \sum_{j=1}^r\left[\gamma_1\left(\dot{x}_i-\dot{x}_{i-j}\right)+\gamma_0\left(x_i-x_{i-j}\right)\right],
\end{align*}
resulting in the propagation coefficient
\begin{equation} \label{eq:Phi_second}
\Phi_2(s)=\frac{a Q_2(s)}{s^2+a r Q_2(s)}=\frac{a\left(\gamma_1 s+\gamma_0\right)}{s^2+a r\left(\gamma_1 s+\gamma_0\right)}.
\end{equation}
Its asymptotic behavior satisfies
\begin{equation}
\left|\Phi_2(0)\right|=\frac{1}{r}, \quad\left|\Phi_2(j \omega)\right| \sim \frac{a \gamma_1}{|\omega|} \rightarrow 0 \text { as } \omega \rightarrow \infty .
\end{equation}
Thus, low-frequency disturbances are attenuated whenever $r\ge2$, matching the structural property observed for the first-order protocol. However, the presence of the velocity-difference term $\gamma_1 s$ increases the effective damping $ar\gamma_1$, introducing additional dynamic shaping that can enhance damping under appropriate gain selection and communication richness.
Consequently, for suitably chosen gains $(\gamma_0,\gamma_1)$, the second-order protocol can satisfy
\begin{equation}
\sup_{\omega\ge0}\left|\Phi_2(j\omega)\right|<1,
\end{equation}
whenever $r\ge2$, with improved attenuation in the mid- and high-frequency ranges relative to the first-order case.

\subsection{String stability of $m^{\text{th}}$-order consensus}
For the generalized $m$-fold integrator, a unified interpretation emerges by examining $\Phi_m(j \omega)$ across frequency ranges. At low frequencies,
\begin{align*}
\left|\Phi_m(0)\right|=\frac{1}{r},
\end{align*}
showing that communication richness alone determines whether disturbances attenuate or persist. Regardless of the protocol order or the gains $\gamma_k$, no consensus law can alter the structural low-frequency limitation imposed by insufficient communication, although higher-order dynamics can significantly influence transient behavior. Thus, under the adopted propagation-based criterion, disturbance attenuation is achievable only if $r\ge2$, and can be attained for suitable gain selections whenever $r\ge2$.

At high frequencies, the leading terms of $s^m$ dominate, and the magnitude of $\Phi_m(j \omega)$ decays to zero at high frequencies, with asymptotic behavior governed by the dominant polynomial terms, which implies strong attenuation of sufficiently high-frequency disturbances for all protocol orders. Increasing the consensus order $m$ primarily influences the transient and frequency-domain behavior, including damping, oscillations, and mid- to high-frequency attenuation. This explains why higher-order consensus protocols can produce smoother transients and improved rejection of oscillatory disturbances, despite sharing the same low-frequency attenuation limit as lower-order protocols.

The unified analysis thus establishes a clear separation between structural and dynamic attenuation. Structural attenuation is governed entirely by the communication richness $r$. Dynamic attenuation, which affects transient smoothness and robustness to higher-frequency disturbances, is governed by the consensus order $m$ and the choice of gains $\gamma_k$. This decomposition generalizes and clarifies several empirical observations in the cooperative driving and consensus literature, providing a coherent theoretical basis for the co-design of communication architectures and distributed control protocols. Moreover, the effectiveness of higher-order consensus in improving transient performance depends on sufficient communication richness, as limited predecessor information may lead to oscillatory or poorly damped responses. These analytical predictions will be validated through time-domain simulations, illustrating how communication richness and consensus order jointly influence disturbance propagation and transient behavior.

\section{Simulation and Results} \label{sec:example}
This section presents numerical simulations illustrating the unified string-stability properties derived in Section~\ref{sec:string_stability}. The objective is not to validate a specific physical vehicle model, but to demonstrate how the consensus order $m$ and communication richness $r$ jointly shape disturbance propagation in an $r$-predecessor leader--follower formation. Both frequency-domain propagation characteristics and time-domain spacing-error responses are examined to show how the theoretical results manifest in practice.

We consider a formation of $n=20$ follower agents, indexed $1: 20$, and a leader whose motion introduces a disturbance into the network. Each follower is simulated using the same $m$th-order integrator model assumed in the analysis, with $m=1,2,3$ depending on the protocol under study and consistent with the closed-loop error dynamics used in Section~\ref{sec:string_stability}. Three consensus protocols are evaluated: first-order $(m=1)$, second-order $(m=2)$, and third-order ( $m=3$ ) consensus. For each order, the communication richness is varied among $r= 1,2,3$, allowing direct comparison across the full $(m, r)$ grid.

\subsection{Simulation Setup}
Each follower satisfies~\eqref{equ:agent_dynamics_morder} and the consensus-based input is implemented according to the general form
\begin{align*}
u_i(t)=a \sum_{k=0}^{m-1} \sum_{j=1}^r\left[\gamma_k\left(x_{i-j}^{(k)}-x_i^{(k)}\right)\right],
\end{align*}
with only the gains corresponding to the selected consensus order $m$ being active. This ensures exact consistency with the closed-loop error dynamics of Section~\ref{sec:string_stability} without introducing artificial transfer-function multiplication or additional filtering.

The simulation is performed directly on the closed-loop error dynamics in Section~\ref{sec:model}, with a unit impulse disturbance $w(t)=\delta(t)$ applied to the highest-order error equation. This provides a common benchmark excitation for comparing propagation behavior across protocol orders. For each consensus order $m$, the disturbance entering the error dynamics is interpreted through the leader’s $m^{\text{th}}$ derivative, consistent with the analytical model. The chosen leader profile therefore provides a common benchmark excitation for comparing propagation behavior across protocol orders.

The gains are deliberately kept identical across all $n$ and across all $r$ to ensure that observed differences arise from intrinsic structural effects. The gains are chosen as $\gamma_0=0.2$, $\gamma_1=1.5$, and $\gamma_2=1.0$ for the third-order protocol, with the appropriate subset of gains applied for the lower-order cases. The coupling gain is set to $a=1.0$. Internal stability is preserved in all simulations, and the structural and dynamic differences are clearly observable. The spacing error is defined as 
\begin{align*}
e_i(t)=x_{i-1}(t)-x_i(t), \quad i=1, \ldots, 20,
\end{align*}
consistent with consensus-based formation control (zero desired spacing). 
For each $(m,r)$, we compute (i) the magnitude of the propagation coefficient $|\Phi_m(j\omega)|$, which serves as the central propagation-based indicator of downstream disturbance attenuation in the present framework, and (ii) the time-domain evolution of spacing errors $e_i(t)$ to visualize how disturbances decay or amplify along the formation.

Although the present study focuses on nominal dynamics, the derived low-frequency limit $\left|\Phi_m(0)\right|=1/r$ is structural and is therefore expected to remain informative under moderate parametric perturbations, provided internal stability is preserved. Quantitative robustness margins under delays, uncertainty, and actuator nonlinearities remain an important topic for future investigation.

\subsection{Results and Discussion}

\begin{figure*}[ht]
\centerline{\includegraphics[width = 1\textwidth]{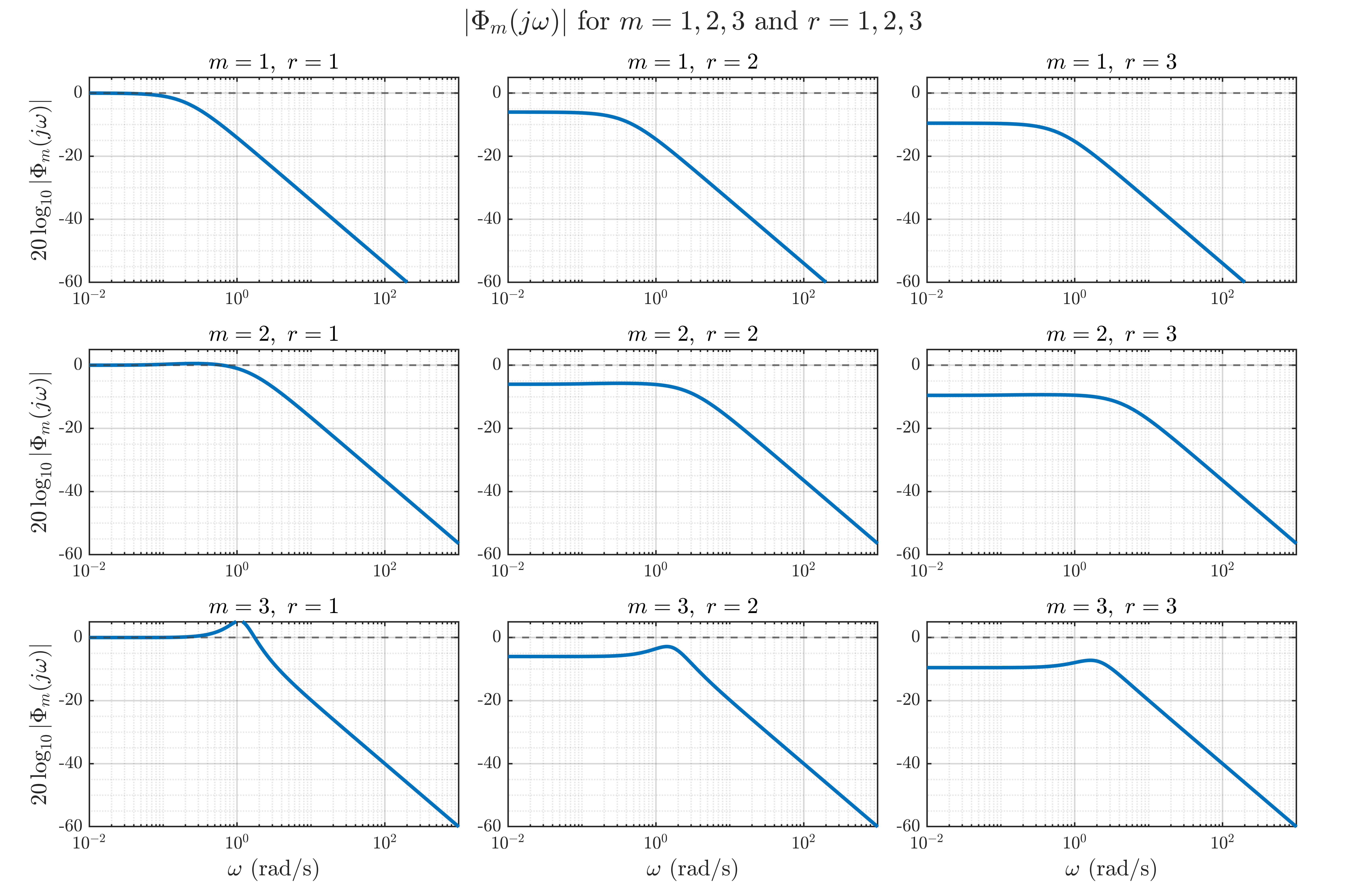}}
\caption{Magnitude of disturbance-propagation coefficient $\left|\Phi_m(j \omega)\right|$ for consensus orders $m=1,2,3$ and communication richness $r=1,2,3$.}
\label{fig:magnitude_plot}
\end{figure*}

Figure~\ref{fig:magnitude_plot} shows the magnitude of the propagation coefficient $\left|\Phi_m(j \omega)\right|$ for all consensus orders and values of $r$. The low-frequency behavior directly corroborates the theoretical conclusion that
\begin{align*}
\left|\Phi_m(0)\right|=\frac{1}{r} \quad \text { for all } \quad m,
\end{align*}
demonstrating that the fundamental string-stability limit is structural and governed solely by communication richness. For $r=1$, all protocols satisfy $\left|\Phi_m(0)\right|=1$, indicating lack of attenuation and amplification of disturbances and the impossibility of low-frequency attenuation; for $r \geq 2$, $\left|\Phi_m(0)\right|<1$, guaranteeing decay of low-frequency disturbances regardless of consensus order.The mid-frequency region reveals that increasing $m$ modifies the spectral shape and can introduce resonant peaks, whose magnitude depends on the communication richness $r$ and gain selection. Increasing $m$ provides additional dynamic shaping, which can improve damping when sufficient communication richness ($r \ge 2$) is available. For sufficiently rich communication, increasing $m$ can improve damping and smoothness; however, when $r=1$, higher-order dynamics may introduce resonant amplification. High-frequency behavior aligns with the analytical expression $\Phi_m(j \omega) \rightarrow 0$ as $\omega \rightarrow \infty$, with higher-order consensus exhibiting steeper decay. These results confirm that the protocol order affects only the mid- and high-frequency attenuation, not the structural low-frequency gain.

\begin{figure}[h]
\centerline{\includegraphics[width = 0.65\textwidth]{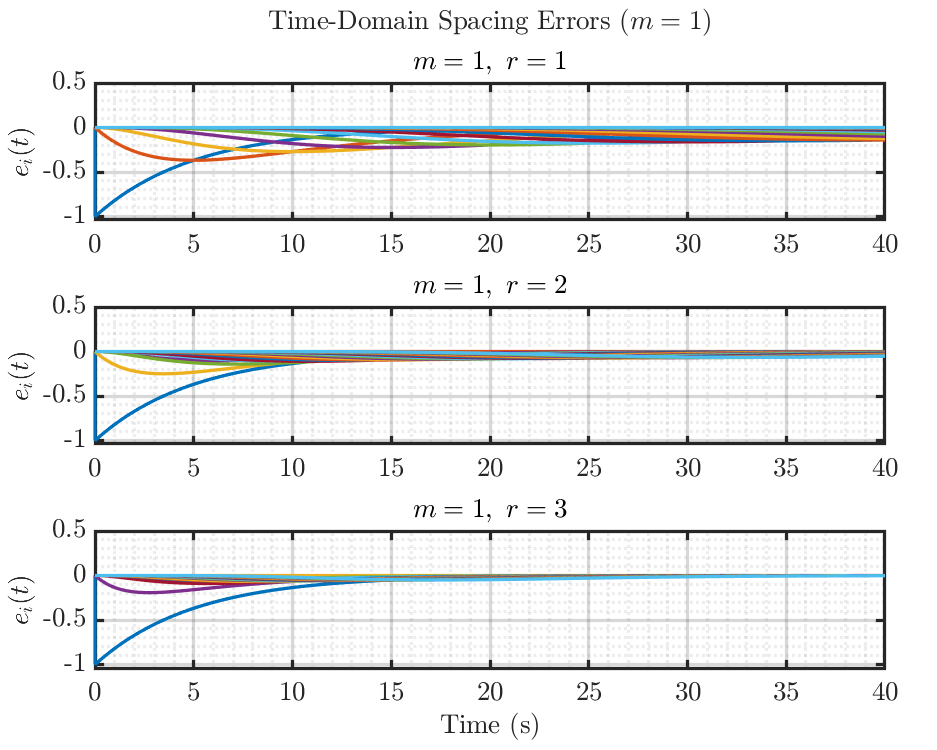}}
\caption{Spacing error trajectories for first-order consensus ( $m=1$ ) under $r=1,2,3$.}
\label{fig:SE_m1}
\end{figure}

\begin{figure}[h]
\centerline{\includegraphics[width = 0.65\textwidth]{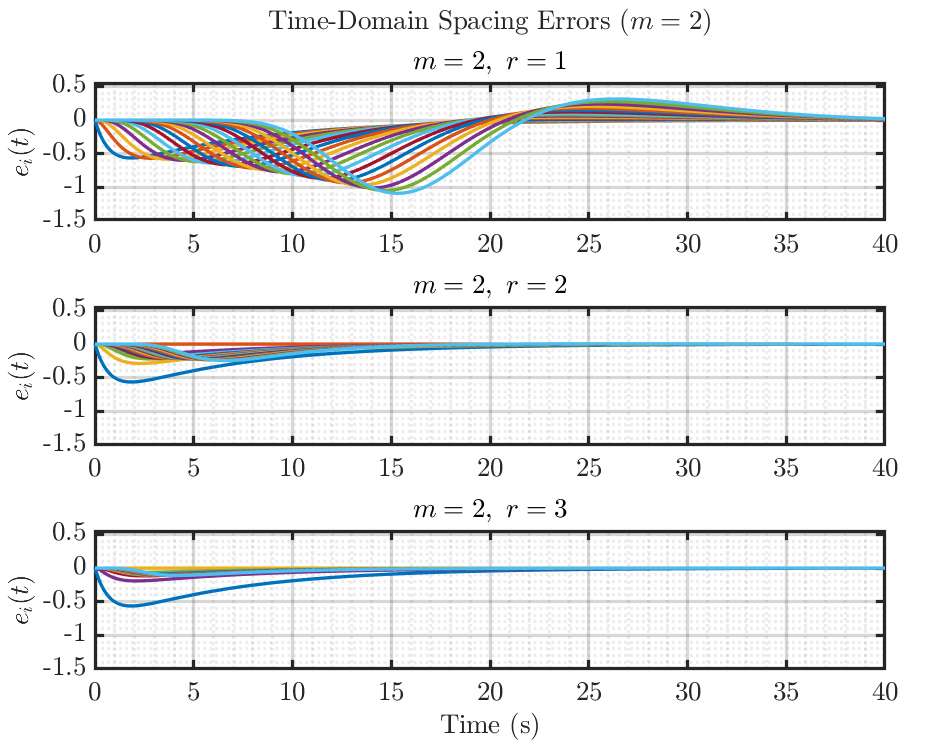}}
\caption{Spacing error trajectories for second-order consensus ( $m=2$ ) under $r=1,2,3$.}
\label{fig:SE_m2}
\end{figure}

\begin{figure}[h]
\centerline{\includegraphics[width = 0.65\textwidth]{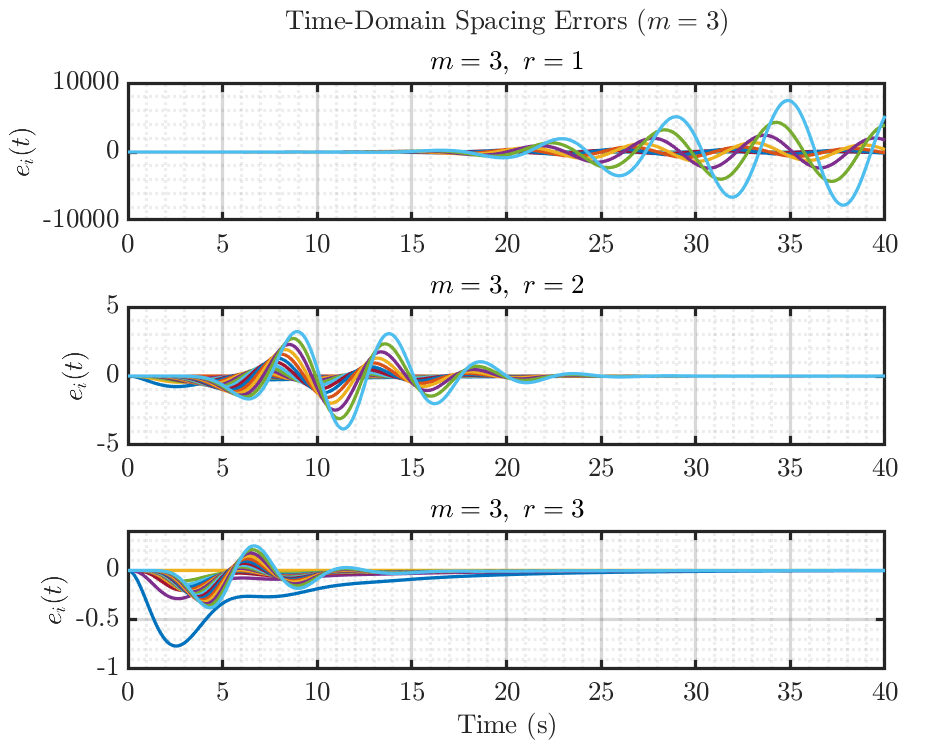}}
\caption{Spacing error trajectories for third-order consensus ( $m=3$ ) under $r=1,2,3$.}
\label{fig:SE_m3}
\end{figure}

Figures~\ref{fig:SE_m1}-\ref{fig:SE_m3} depict the spacing-error evolution across followers for $m=1,2,3$. The time responses exhibit the same qualitative phenomena predicted by the propagation coefficient.

For first-order consensus, the formation with $r=1$ exhibits non-decaying spacing errors, consistent with the lack of attenuation under single-predecessor communication. Increasing the communication richness to $r=2$ or $r=3$ produces clear downstream attenuation, though the first-order consensus exhibits monotonic but relatively slow decay due to limited dynamic shaping. This confirms that first-order consensus can be string stable only when $r \geq 2$.

For second-order consensus, the introduction of velocity coupling substantially suppresses oscillations. Although $r=1$ again prevents full attenuation, the responses are smoother and settle faster compared to first-order consensus. When $r \geq 2$, spacing errors decay rapidly, with significantly improved transient performance relative to the first-order case. This matches the theoretical prediction that increasing $m$ affects mid-frequency dynamics.

For third-order consensus, the influence of higher-order dynamics becomes more pronounced. When $r=1$, the responses exhibit significant oscillatory amplification, indicating that increasing the consensus order without sufficient communication richness can degrade performance. In contrast, for $r \geq 2$, the formation exhibits rapid attenuation and improved transient behavior. This demonstrates that higher-order consensus provides enhanced dynamic shaping only when adequate communication information is available.

Across all simulations, the central theoretical conclusions of Section~\ref{sec:string_stability} are
\begin{enumerate}
    \item Communication richness $r$ solely determines string stability: The low-frequency gain $1 / r$ is visible in both the frequency responses and in the spacing-error propagation. When $r=1$, disturbances do not attenuate; when $r \geq 2$, they decay along the formation.
    \item Consensus order $m$ shapes transient quality but not string stability: Higher-order protocols provide additional dynamic shaping, but improved damping and smoothness depend on sufficient communication richness and appropriate gain selection. However, no choice of consensus order or gains can alter the low-frequency structural limit imposed by $r$.
    \item The qualitative trends observed across all simulations precisely reflect the magnitude characteristics of $\Phi_m(j \omega)$, confirming that the entire disturbance-propagation mechanism is captured by this single scalar function.
\end{enumerate}

Together, the simulations reinforce the main insight of the paper: string stability in consensus-driven formations is fundamentally a topological property, not a controller-order property, and higher-order consensus improves performance only in dynamic, not structural, ways. Notably, the simulations also reveal that increasing the consensus order without sufficient communication richness can degrade performance, highlighting the need for co-design of communication topology and protocol dynamics.

\section{Conclusion}
This paper presented a unified frequency-domain framework for analyzing string stability of leader--follower formations governed by first-, second-, and $m^{\text{th}}$ order consensus protocols under an $r$-predecessor communication topology. Using a common propagation coefficient $\Phi_m(s)$, we established that all consensus orders share the same structural low-frequency gain $\left|\Phi_m(0)\right|=1/r$, leading to the key finding that downstream disturbance attenuation is structurally achievable if and only if $r\ge2$. Thus, communication richness fundamentally determines disturbance attenuation, whereas increasing the consensus order refines transient and mid-frequency behavior but cannot alter the low-frequency limitation.

Simulation results across $m=1,2,3$ and $r=1,2,3$ confirmed these analytical insights: higher-order consensus can improve damping and smoothness when supported by sufficient predecessor information, while only $r$ controls whether disturbances amplify or decay along the formation. The simulations further show that increasing protocol order without sufficient communication richness can lead to oscillatory amplification, highlighting the need for joint topology--protocol design. Together, the theory and simulations demonstrate a clear separation between structural attenuation dictated by topology and dynamic attenuation shaped by protocol order. 

This work considered a structured $r$-predecessor topology, which permits an exact scalar propagation characterization. For more general directed graphs, disturbance propagation depends on the full Laplacian spectrum and potentially heterogeneous modal interactions, making explicit string-stability criteria substantially more challenging. Extending the framework to more general directed graph topologies, heterogeneous agent dynamics, communication delays, and adversarial disturbances remains a natural direction for future research \cite{10644849,ZHENG2023105552}.

\bibliographystyle{unsrtnat}
\bibliography{references}

@article{ZHENG2023105552,
title = {Moving-horizon false data injection attack design against cyber–physical systems},
journal = {Control Engineering Practice},
volume = {136},
pages = {105552},
year = {2023},
author = {Yu Zheng and Sridhar Babu Mudhangulla and Olugbenga Moses Anubi}
}

@INPROCEEDINGS{10644849,
  author={Rajarajan, Naveen Kumar and Mudhangulla, Sridhar Babu and Anubi, Olugbenga Moses},
  booktitle={2024 American Control Conference (ACC)}, 
  title={Passive Stability and Adaptive Control of Teleoperated System Using Wave Variables and Predictor Techniques}, 
  year={2024},
  pages={4365-4371}
  }

@article{feng2019string,
  title={String stability for vehicular platoon control: Definitions and analysis methods},
  author={Feng, Shuo and Zhang, Yi and Li, Shengbo Eben and Cao, Zhong and Liu, Henry X and Li, Li},
  journal={Annual Reviews in Control},
  volume={47},
  pages={81--97},
  year={2019},
  publisher={Elsevier}
}

@article{ren2008consensus,
  title={Consensus algorithms for double-integrator dynamics},
  author={Ren, Wei and Beard, Randal W},
  journal={Distributed Consensus in Multi-vehicle Cooperative Control: Theory and Applications},
  pages={77--104},
  year={2008},
  publisher={Springer}
}

@article{bian2019reducing,
  title={Reducing time headway for platooning of connected vehicles via V2V communication},
  author={Bian, Yougang and Zheng, Yang and Ren, Wei and Li, Shengbo Eben and Wang, Jianqiang and Li, Keqiang},
  journal={Transportation Research Part C: Emerging Technologies},
  volume={102},
  pages={87--105},
  year={2019},
  publisher={Elsevier}
}

@inproceedings{ren2006high,
  title={High-order consensus algorithms in cooperative vehicle systems},
  author={Ren, Wei and Moore, Kevin and Chen, YangQuan},
  booktitle={2006 IEEE International Conference on Networking, Sensing and Control},
  pages={457--462},
  year={2006},
  organization={IEEE}
}

@article{wieland2008consensus,
  title={On consensus in multi-agent systems with linear high-order agents},
  author={Wieland, Peter and Kim, Jung-Su and Scheu, Holger and Allg{\"o}wer, Frank},
  journal={IFAC Proceedings Volumes},
  volume={41},
  number={2},
  pages={1541--1546},
  year={2008},
  publisher={Elsevier}
}

@article{besselink2017string,
  title={String stability and a delay-based spacing policy for vehicle platoons subject to disturbances},
  author={Besselink, Bart and Johansson, Karl H},
  journal={IEEE Transactions on Automatic Control},
  volume={62},
  number={9},
  pages={4376--4391},
  year={2017},
  publisher={IEEE}
}

@article{pare2019networked,
  title={Networked model for cooperative adaptive cruise control},
  author={Par{\'e}, Philip E and Hashemi, Ehsan and Stern, Raphael and Sandberg, Henrik and Johansson, Karl Henrik},
  journal={IFAC-PapersOnLine},
  volume={52},
  number={20},
  pages={151--156},
  year={2019},
  publisher={Elsevier}
}

@inproceedings{barooah2005error,
  title={Error amplification and disturbance propagation in vehicle strings with decentralized linear control},
  author={Barooah, Prabir and Hespanha, Joao P},
  booktitle={Proceedings of the 44th IEEE conference on Decision and Control},
  pages={4964--4969},
  year={2005},
  organization={IEEE}
}

@article{ploeg2013lp,
  title={Lp string stability of cascaded systems: Application to vehicle platooning},
  author={Ploeg, Jeroen and Van De Wouw, Nathan and Nijmeijer, Henk},
  journal={IEEE Transactions on Control Systems Technology},
  volume={22},
  number={2},
  pages={786--793},
  year={2013},
  publisher={IEEE}
}

@article{swaroop2002string,
  title={String stability of interconnected systems},
  author={Swaroop, Darbha and Hedrick, J Karl},
  journal={IEEE transactions on automatic control},
  volume={41},
  number={3},
  pages={349--357},
  year={2002},
  publisher={IEEE}
}

@article{ren2007information,
  title={Information consensus in multivehicle cooperative control},
  author={Ren, Wei and Beard, Randal W and Atkins, Ella M},
  journal={IEEE Control systems magazine},
  volume={27},
  number={2},
  pages={71--82},
  year={2007},
  publisher={IEEE}
}

@article{olfati2004consensus,
  title={Consensus problems in networks of agents with switching topology and time-delays},
  author={Olfati-Saber, Reza and Murray, Richard M},
  journal={IEEE Transactions on automatic control},
  volume={49},
  number={9},
  pages={1520--1533},
  year={2004},
  publisher={IEEE}
}

@article{naus2010string,
  title={String-stable CACC design and experimental validation: A frequency-domain approach},
  author={Naus, Gerrit JL and Vugts, Rene PA and Ploeg, Jeroen and van De Molengraft, Marinus JG and Steinbuch, Maarten},
  journal={IEEE Transactions on vehicular technology},
  volume={59},
  number={9},
  pages={4268--4279},
  year={2010},
  publisher={IEEE}
}

@book{godsil2013algebraic,
  title={Algebraic graph theory},
  author={Godsil, Chris and Royle, Gordon F},
  volume={207},
  year={2013},
  publisher={Springer Science \& Business Media}
}

@inproceedings{zegers2016consensus,
  title={Consensus-based bi-directional CACC for vehicular platooning},
  author={Zegers, Jeroen C and Semsar-Kazerooni, Elham and Ploeg, Jeroen and van de Wouw, Nathan and Nijmeijer, Henk},
  booktitle={2016 American Control Conference (ACC)},
  pages={2578--2584},
  year={2016},
  organization={IEEE}
}

@article{overview_salek_2025,
	title = {An Overview of Automated Vehicle Longitudinal Platoon Formation Strategies},
	doi = {10.1145/3736642},
	author = {Salek, M Sabbir and Thakur, Mugdha Basu and Krishna, Pardha Sai and Chowdhury, Mashrur and Schmid, Matthias and Murray‐Tuite, Pamela and Khan, Sakib Mahmud and Krovi, Venkat},
	journal = {ACM Journal on Autonomous Transportation Systems},
	year = {2025},
	litmapsId = {287121697}
}

@article{performance_hansson_2024,
	title = {Performance Bounds for Multi-Vehicle Networks With Local Integrators},
	doi = {10.1109/lcsys.2024.3518397},
	author = {Hansson, Jonas and Tegling, Emma},
	journal = {IEEE Control Systems Letters},
	year = {2024},
	litmapsId = {282484109}
}

@article{vehicular_stdli_2017,
	title = {Vehicular Platoons in cyclic interconnections with constant inter-vehicle spacing},
	doi = {10.1016/j.ifacol.2017.08.449},
	author = {Stüdli, Sonja and Seron, Maria M. and Middleton, Richard H.},
	journal = {IFAC-PapersOnLine},
	year = {2017},
	litmapsId = {21514547}
}

@article{overview_knorn_2016,
	title = {Overview: Collective Control of Multiagent Systems},
	doi = {10.1109/tcns.2015.2468991},
	author = {Knorn, S. and Chen, Zhiyong and Middleton, R.},
	journal = {IEEE Transactions on Control of Network Systems},
	year = {2016},
	litmapsId = {75607802}
}

@article{disturbance_seiler_2004,
	title = {Disturbance propagation in vehicle strings},
	doi = {10.1109/tac.2004.835586},
	author = {Seiler, P. and Pant, A. and Hedrick, J.},
	journal = {IEEE Transactions on Automatic Control},
	year = {2004},
	litmapsId = {23001964}
}

@article{output_zhu_2025,
	title = {Output consensus of heterogeneous switched nonlinear multi-agent systems with sampled-data communication},
	doi = {10.1108/ria-01-2025-0009},
	author = {Zhu, Jingyi and Zou, Wencheng},
	journal = {Robotic Intelligence and Automation},
	year = {2025},
	litmapsId = {288369012}
}

@article{fully_mazouchi_2021,
	title = {Fully Heterogeneous Containment Control of a Network of Leader–Follower Systems},
	doi = {10.1109/tac.2021.3130878},
	author = {Mazouchi, Majid and Tatari, Farzaneh and Kiumarsi-Khomartash, Bahare and Modares, H.},
	journal = {IEEE Transactions on Automatic Control},
	year = {2021},
	litmapsId = {254524067}
}

@article{constrained_zhou_2019,
	title = {Constrained Consensus in Continuous-Time Multi-Agent Systems},
	author = {Zhou, Zheqing},
	year = {2019},
	litmapsId = {70467421}
}

@article{multiagent_jiang_2022,
	title = {Multi-agent consensus with heterogeneous time-varying input and communication delays in digraphs},
	doi = {10.1016/j.automatica.2021.109950},
	author = {Jiang, Wei and Liu, Kun and Charalambous, Themistoklis},
	journal = {at - Automatisierungstechnik},
	year = {2022},
	litmapsId = {98276398}
}

@article{distributed_huang_2020,
	title = {Distributed output feedback consensus control of networked homogeneous systems with large unknown actuator and sensor delays},
	doi = {10.1016/j.automatica.2020.109249},
	author = {Huang, Ran and Ding, Zhengtao and Cao, Zhengcai},
	journal = {at - Automatisierungstechnik},
	year = {2020},
	litmapsId = {46360580}
}

@book{zhou1998essentials,
  title={Essentials of robust control},
  author={Zhou, Kemin and Doyle, John Comstock},
  volume={104},
  year={1998},
  publisher={Prentice hall Upper Saddle River, NJ}
}

\end{document}